\documentstyle[12pt,twoside]{article}

\newcommand{\ve}{\vspace{3mm}}

\font\ninerm=cmr9

\font\nineit=cmti9
\font\ninebf=cmbx9

\newcommand{\be}{\begin{equation}}
\newcommand{\ee}{\end{equation}}
\newcommand{\bea}{\begin{eqnarray}}
\newcommand{\eea}{\end{eqnarray}}
\newcommand{\bean}{\begin{eqnarray*}}
\newcommand{\eean}{\end{eqnarray*}}

\newfont\frak{eufm10}

\ve
\newcommand{\R}{I\!\! R}
\newcommand{\mN}{I\!\!N}

\newcommand{\mQ}{\mbox{}\; l\!\!\!Q}

\newcommand{\mC}{\mbox{}\; l\!\!\!C}
\newcommand{\mR}{I\!\!R}
\newcommand{\spar}{\vskip0.3cm}

\newtheorem{proposition}{Proposition}
\newtheorem{theorem}[proposition]{Theorem}

\newtheorem{lemma}[proposition]{Lemma}
\newtheorem{definition}[proposition]{Definition}
\newtheorem{remark}[proposition]{Remark}

\begin{document}


\markboth{Real equation solving}{Real equation solving}

\title{\bf Polar Varieties, Real Equation Solving and Data-Structures: The 
 hypersurface case \footnote{Research partially supported by the Spanish 
 government grant DGICT, PB93-0472-C02-02}}

\author{\sc B. Bank $^{1}$, M. Giusti $^{2}$, J. Heintz $^{3}$,\\[.5cm] \sc
                  G. M. Mbakop $^{1}$}

\maketitle
\begin{center} {\em Dedicated to Shmuel Winograd} \end{center}

\addtocounter{footnote}{1}\footnotetext{Humboldt-Universit\"at zu Berlin,
Untern den Linden 6, D--10099 Berlin, Germany.\\
\hbox{\hspace{.5cm}}
{bank@mathematik.hu-berlin.de, mbakop@mathematik.hu-berlin.de} }

\addtocounter{footnote}{1}\footnotetext{GAGE, Centre de
Math\'ematiques, \'Ecole Polytechnique, F-91228 Palaiseau Cedex, France.\\
\hbox{\hspace{.5cm}} {giusti@ariana.polytechnique.fr} }

\addtocounter{footnote}{1}\footnotetext{ Dept. de Matem\'aticas,
Estad\'{\i}stica y Computaci\'on, Facultad de Ciencias, Universidad de\\
\hbox{\hspace{.5cm}} Cantabria, E-39071 Santander, Spain. 
{heintz@matsun1.unican.es} }



\ve
\begin{abstract}

In this paper we apply for the first time a new method for
multivariate equation solving which was developed in \cite{gh1},
\cite{gh2}, \cite{gh3} for complex root determination  to the {\em real\/} case. 
Our main result
concerns the problem of finding at least one representative point for each
connected component of a real compact and smooth hypersurface. \spar

The basic algorithm of \cite{gh1}, \cite{gh2}, \cite{gh3} yields
a new method for symbolically solving zero-dimensional
polynomial equation systems over the complex numbers. One
feature of central importance of this algorithm is the use of a
problem--adapted data type represented by the data structures
arithmetic network and straight-line program (arithmetic
circuit).  The algorithm finds the complex solutions of any
affine zero-dimensional equation system in non-uniform
sequential time that is {\em polynomial\/} in the length of the
input (given in straight--line program representation) and an
adequately defined {\em geometric degree of the equation
system}. \spar Replacing the notion of geometric degree of the
given polynomial equation system by a suitably defined {\em real
(or complex) degree\/} of certain polar varieties associated to
the input equation of the real hypersurface under consideration,
we are able to find for each connected component of the
hypersurface a representative point (this point will be given in
a suitable encoding). The input equation is supposed to be given
by a straight-line program and the (sequential time) complexity
of the algorithm is polynomial in the input length and the
degree of the polar varieties mentioned above.
\spar


\end{abstract}


\begin{section}{Introduction}

The present article is strongly related to the main complexity
results and algorithms in \cite{gh1}, \cite{gh2}, \cite{gh3}.
Whereas the algorithms developed in these papers concern solving
of polynomial equation systems over the complex numbers, here we
deal with the problem of real solving. More precisely, we
consider the particular problem of finding real solutions of a
single equation $f(x)=0$, where $f$ is an $n$--variate
polynomial of degree $d \ge 2$ over the rationals which is
supposed to be a regular equation of a compact and smooth
hypersurface of $\mR^n$. Best known complexity bounds for this
problem over the reals are of the form $d^{O(n)}$, counting
arithmetic operations in $\mQ$ at unit cost (see
\cite{grivo1}, \cite{Griegorev87}, \cite{sole}, \cite{hroy}, \cite{hroy1},
\cite{hroy2}, \cite{canny}, \cite{Renegar1}, \cite{rene}, \cite{basu}).  \spar

Complex root finding methods cannot be applied directly to real
polynomial equation solving just by looking at the complex
interpretation of the input system. If we want to use a complex
root finding method for a problem over the reals, some previous
adaptation or preprocessing of the input data becomes necessary.
In this paper we show that certain {\em polar varieties\/}
associated to our input affine hypersurface possess specific
geometric properties, which permits us to adapt the complex main
algorithm designed in the papers \cite{gh1}, \cite{gh2},
\cite{gh3} to the real case.
\spar

This algorithm is of {\em intrinsic type}, which means that it
allows to distinguish between semantical and syntactical
properties of the input system in order to profit from both for
an improvement of the complexity estimates compared with more
"classical" procedures (as e.g. \cite{her},
\cite{seid}, \cite{buch}, \cite{hew}, \cite{he}, \cite{laz1}, \cite{laz2},
\cite{ChGr}, \cite{can-gal-he}, \cite{dickeal}, \cite{canny},
\cite{giu-he92}, \cite{krick-pardo:CRAS}, \cite{krick-pardo1},
\cite{Chistov95}). The papers \cite{gh1}, \cite{gh2} and \cite{gh3} show that
the {\em geometric degree of the input system\/} is associated
with the intrinsic complexity of solving the system
algorithmically when the complexity is measured in terms of the
number of arithmetic operations in $\mQ$.  The paper \cite{gh1}
is based on the somewhat unrealistic complexity model in which
certain {\em FOR\/} instructions executable in parallel count at
unit cost. This drawback of the complexity model is corrected in
the paper
\cite{gh2} at the price of introducing algebraic parameters in the
straight--line programs and arithmetic networks occurring there.
These algebraic parameters are finally eliminated in the paper
\cite{gh3}, which contains a procedure satisfying our complexity
requirement and is completely rational.  \spar

We show that the algorithmic method of the papers \cite{gh1},
\cite{gh2} and
\cite{gh3} is also applicable to the problem of (real) root finding in the
case of a compact and smooth hypersurface of $\mR^n$, given by
an $n-$variate polynomial $f$ of degree $d$ with rational
coefficients which represents a regular equation of that
hypersurface. It is possible to design an algorithm of {\em
intrinsic type\/} using the same data structures as in
\cite{gh3}, namely arithmetic networks and straight-line programs over $\mQ$
(the straight--line programs -- which are supposed to be
division--free -- are used for the coding of input system,
intermediate results and output).  In the complexity estimates
the notion of {\em (geometric) degree of the input system\/} of
\cite{gh1}, \cite{gh2}, \cite{gh3} has then to be replaced by
the {\em (complex or real) degree of the polar varieties} which
are associated to the input equation.
\spar 

The basic computation model used in our algorithm will be that
of an arithmetic network with parameters in $\mQ$ (compare with
\cite{gh3}).  Our first complexity result is the following:
\spar

{\em There is an arithmetic network of size $(nd\delta
L)^{O(1)}$ with parameters in the field of the rational numbers
which finds at least one representative point in every connected
component of a smooth compact hypersurface of $\mR^n$ given by a
regular equation $f \in \mQ[X_1, \ldots , X_n]$ of degree $d \ge
2$. Here $L$ denotes the size of a suitable straight-line
program which represents the input of our procedure coding the
input polynomial $f$. Moreover, $\delta$ denotes the
maximal geometric degree of suitably defined polar varieties
associated to the input equation $f$.}
\spar

Our second complexity result relies on two algorithmic
assumptions which are very strong in theory, but hopefully not
so restrictive in practice. We assume now that a factorization
procedure for univariate polynomials over $\mQ$ being
"polynomial" in a suitable sense (e.g. counting arithmetic
operations in $\mQ$ at unit cost) is available and that we are
able (also at polynomial cost) to localize regions where a given
multivariate polynomial has "many" real zeros (if there exist
such regions). This second assumption may be replaced by the
following more theoretical one (which, however, is simpler to
formulate precisely): we suppose that we are able to decide in
polynomial time whether a given multivariate polynomial has a
real zero (however we do {\em not\/} suppose that we are able to
exhibit such a zero if there exists one).  We call an arithmetic
network {\em extended\/} if it uses subroutines of these two
types.
\spar

{\em Let notations and assumptions be as before. Suppose
furthermore that $f$ represents a regular equation of a
non-empty smooth and compact real hypersurface. Then there
exists an extended arithmetic network which finds at least one
representative point for each connected component of the real
hypersurface given by $f$. The size of this arithmetic network
is $(nd\delta^* L)^{O(1)}$, where $\delta^*$ denotes the
suitably defined maximal real degree of the polar varieties
mentioned above.}
\spar

Complexity results in a similar sense for the specific problem
of {\em numerical\/} polynomial equation solving can be found in
\cite{ShSm1}, following an approach initiated in \cite{ShSm93a},
\cite{ShSm93b},
\cite{ShSm93c}, \cite{ShSm93d} (see also \cite{Dedieu1}, \cite{Dedieu2}).
In the same sense one might also want to mention \cite{CaEm95}
and
\cite{Emiris} as representative contributions for the sparse viewpoint.
For more details we refer the reader to \cite{pardo} and
\cite{gh3} and the references cited therein.

\end{section} 

\begin{section}{Polar Varieties}

As usual, let $\mQ, \mR$ and $\mC$ denote the field of rational,
real and complex numbers, respectively. The affine n--spaces
over these fields are denoted by $\mQ^n, \mR^n$ and $\mC^n$,
respectively. Further, let $\mC^n$ be endowed with the Zariski
topology of $\mQ-$definable algebraic sets, where a closed set
consists of all common zeros of a finite number of polynomials
with coefficients in $\mQ$. Let $W \subset {\mC}^n$ be a closed
subset with respect to this topology and let $W= C_1\cup\cdots
\cup C_s$ be its decomposition into irreducible components with
respect to the same topology. Thus $W, C_1,\ldots,C_s$ are
algebraic subsets of~ ${\mC}^n$.  We call $W$ equidimensional if
all its irreducible components $C_1, \ldots, C_s$ have the same
dimension.
\spar

In the following we need the notion of (geometric) degree of an
affine algebraic variety. Let $W$ be an equidimensional Zariski
closed subset of $\mC^n$. If $W$ is zero--dimensional, the {\em
degree\/} of $W$, denoted by $deg W$, is defined as the
cardinality of $W$ (neither multiplicities nor points at
infinity are counted). If $W$ is of positive dimension $r$, then
we consider the collection ${\cal M}$ of all $(n-r)$-dimensional
affine linear subspaces, given as the solution set in $\mC^n$ of
a linear equation system $L_1 = 0, \ldots, L_{r} = 0$ where for
$1\leq k\leq r$ the equation $L_k$ is of the form $L_{k} =
\sum_{j=1}^{n} a_{kj} x_j + a_{k0}$ with $a_{kj}$ being
rational.  Let ${\cal M}_{W}$ be the subcollection of ${\cal M}$
consisting of all affine linear spaces $H \in {\cal M}$ such
that the affine variety $H \cap W$ satisfies $H \cap W \not=
\emptyset$ and $dim(H \cap W) = 0$. Then the geometric degree of
$W$ is defined as $deg W := max\{ deg ( W \cap H ) | H \in {\cal
M}_W \}$.  \spar

For an {\em arbitrary\/} Zariski closed subset $W$ of $\mC^n$,
let $W= C_1\cup\cdots \cup C_s$ be its decomposition into
irreducible components.  As in \cite{he} we define its geometric
degree as $deg W := deg C_1 + \cdots + deg C_s$. Let $W$ be a
Zariski closed subset of $\mC^n$ of dimension $n-i$ given by a
regular sequence of polynomials $f_1, \ldots, f_i \in \mQ[X_1,
\ldots, X_n]$.

\begin{definition}\label{def1}

For $1 \le j \le s$, the irreducible component $C_j$ is called a
{real component\/} of $W$ if the real variety $C_j\cap \mR^n$
contains a smooth point of $C_j$. Let us write

\[ I := \{ j \in \mN | 1 \le j \le s, \hskip 3pt C_j \hskip 3pt{\hbox {
is a real component of $W$}} \}. \]

Then the (complex) affine variety $W^\ast := \bigcup \limits_{j
\in I} C_j$ is called the  real part of $W$. We call $deg^{\ast}
W := deg W^{\ast} =
\sum\limits_{j \in I} deg C_j$ the real degree of the algebraic set
$W$.

\end{definition}

\begin{remark}

(i) $deg^{\ast} W= 0$ holds if and only if the real part
$W^{\ast}$ of $W$ is empty.\\

(ii)  Note that  "smooth point of $C_j$" in Definition 1 is
somewhat ambiguous and should be interpreted following the
context.  Thus "smooth point of $C_j$" may just mean that the
tangent space of $C_j$ is of dimension $(n-i)$ at such a point,
or, more restrictively, it may mean that the hypersurfaces
defined by the polynomials $f_1, \ldots, f_i$ intersect
transversally in such a  point.
\end{remark}

\begin{proposition}\label{prop3}

Let $f \in \mQ[X_1,\ldots,X_n]$ be a non-constant and
square-free polynomial and let $W := \{ x \in \mC^n | f(x) =
0\}$ be the set of complex zeros of the equation $f(x) =0$.
Furthermore, consider for any fixed $i, 0 \le i <n$, the complex
variety

$$\widetilde{W}_i  :=  \left \{ x \in \mC^n |  f(x) = {{\partial
f(x)} \over {\partial X_1}} = \cdots = {{\partial f(x)} \over
{\partial X_i}} = 0
\right \}$$

(here $\widetilde{W}_0$ is understood to be $W$). Suppose that
the variables $X_1, \ldots, X_n$ are in generic position with
respect to $f$. Then any point of $\widetilde{W}_i$ being a
smooth point of $W$ is also a smooth point of $\widetilde{W}_i$.
More precisely, at any such point the Jacobian of the equation
system $f=\frac{\partial f(x)}{\partial X_1} = \cdots =
\frac{\partial f(x)}{\partial X_i} =0$ has maximal rank, i.e., the
hypersurfaces defined by the polynomials $ f, {{\partial f}
\over {\partial X_1}}, \ldots,{{\partial f} \over {\partial
X_i}}$ intersect transversally in this point.

\end{proposition}
\spar

{\bf Proof:}

Consider the non-singular linear transformation $x
\longleftarrow A^{(i)} y$, where the new variables are $y =
(Y_1, \ldots, Y_n)$. Suppose that $A^{(i)}$ is given in the form

\begin{equation}
\left( \begin{array}{ll}
I_{i,i} & 0_{i,n-i} \\ (a_{kl})_{n-i,i} & I_{n-i,n-i}
\end{array} \right)
\end{equation}

where $I_{i,i}$ and $0_{i, (n-i)}$ denote the $i \times i$  unit
and the $i
\times (n-i)$ zero matrix, respectively, and where $a_{kl}$ are arbitrary
complex numbers for $i+1 \le k \le n$ and $1 \le l \le i$. Since
the square matrix $A^{(i)}$ has full rank, the transformation $x
\longleftarrow A^{(i)} y$ defines a linear change of
coordinates. In the new coordinates, the variety
$\widetilde{W}_i$ takes the form

$$\widetilde{W}_i \!:=\!\! \left \{ \!y \in \mC^n | f(y)\! =\!
{{\partial f(y)} \over {\partial Y_1}}\! +\!\!\!\sum_{j = i+1}^n
\!\!a_{j1} {{\partial f(y)} \over {\partial Y_j}}\! =\! \cdots
\!=\! {{\partial f(y)} \over {\partial Y_i}}\! +\!\!\!\sum_{j =
i+1}^n \!\!a_{ji} {{\partial f(y)} \over {\partial Y_j}}\!=\!
0\!\! \right \}.$$

The coordinate transformation given by $A^{(i)}$ induces a
morphism of affine spaces $\Phi_i  : \mC^n \times \mC^{(n-i) i}
\longrightarrow
\mC^{i+1}$ defined by

$$\Phi_i \left ( Y_1, \ldots, Y_i, \ldots, Y_n, a_{i+1, 1},
\ldots, a_{n,1},
\ldots, a_{i+1,i}, \ldots, a_{n, i} \right ) = $$ $$\left ( f, {{\partial
f} \over {\partial Y_1}} + \sum_{j = i+1}^n a_{j 1} {{\partial
f} \over {\partial Y_j}}, \ldots, {{\partial f} \over {\partial
Y_i}} +\sum_{j = i+1}^n a_{j i} {{\partial f} \over {\partial
Y_j}} \right ).$$

For the moment let
\[
\alpha := (\alpha_1, \ldots, \alpha_{n + (n-i) i} ) := (Y_1, \ldots, Y_n,
a_{i+1,  1}, \ldots, a_{n,i}) \in \mC^n \times \mC^{(n-i) i}
\]
Then the Jacobian matrix $J (\Phi_i) (\alpha )$
of $\Phi_i $ in $ \alpha $ is given by
\[
J (\Phi_i) (\alpha ) =
\left ( \begin{array}{ccllclccl}
 {{\partial f}\over{\partial Y_1}}& \cdots &
 {{\partial f} \over {\partial Y_n}}  & 0 & \cdots & 0 & \cdots & \cdots & 0 \\
 \ast & \cdots&   \ast  & {{\partial f} \over {\partial Y_{i+1}}} & \cdots
& {{\partial f} \over {\partial Y_n}}
 & 0 \cdots  & \vdots & 0\\
 \vdots & & \vdots & \ddots & \ddots & 0 & \cdots & \ddots & 0 \\
 \ast & \cdots& \ast & 0 \cdots  & 0 \cdots & \cdots &{{\partial f} \over
{\partial Y_{i+1}}}
  & \cdots & {{\partial f} \over {\partial Y_n}}
                       \end{array}    \right ) (\alpha )
\]
Suppose that we are given a point $\alpha^0 = (Y_1^0, \ldots,
Y_n^0, a_{i+1, 1}^0, \ldots a_{n, i}^0)$ which belongs to the
fiber $\Phi_i^{-1} (0)$ and suppose that $(Y_1^0, \ldots,
Y_n^0)$ is a point of the hypersurface $W$ in which the equation
$f$ is regular (i.e., we suppose that not all partial
derivatives of $f$ vanish in that point). Let us consider the
Zariski open neighbourhood $\cal U$ of $(Y_1^0, \ldots, Y_n^0)$
consisting of all points of\ $\mC^n$ in which at least one
partial derivative of $f$ does not vanish. We claim now that the
restricted map

\[ \Phi_i : {\cal U} \times \mC^{(n-i)i} \longrightarrow \mC^{i+1} \]

is transversal to the origin $ 0 =(0, \ldots, 0)$ of $\mC^{i+1}$. In order
to prove this assertion we consider an arbitrary point $\alpha = (Y_1, \ldots,
Y_n, a_{i+1, 1}, \ldots, a_{n, i})$ of $\;{\cal U} \times \mC^{(n-i) i}$
which satisfies $\Phi_i (\alpha) = 0$. Thus $(Y_1, \ldots, Y_n)$ belongs to
$\;{\cal U} \cap W$ and is therefore a point of the hypersurface $W$ in which
the equation $f$ is regular. Let us now show that the Jacobian matrix of
$\Phi_i$ has maximal rank in $\alpha$. If this is not the case, the partial
derivatives ${{\partial f} \over {\partial Y_{i+1}}}, \ldots, {{\partial f}
\over {\partial Y_{n}}}$ must vanish in the point $(Y_1, \ldots, Y_n)$. 
Then the relation $\Phi_i (\alpha) = 0$ implies that the derivatives
${{\partial f} \over {\partial Y_1}}, \ldots, {{\partial f} \over {\partial
Y_{i}}}$ at the point $(Y_1, \ldots, Y_n)$ vanish, too.
\spar

This contradicts the fact that the equation $f$ is regular in that point.
Therefore the Jacobian matrix of $\Phi_i$ has maximal rank in $\alpha$, which
means that $\alpha$ is a regular point of $\Phi_i$. Since $\alpha$ was
an arbitrary point of $\Phi_i^{-1}(0) \cap ({\cal U} \times \mC^{(n-i) i})$,
our claim follows. Applying the algebraic--geometric form of the Weak
Transversality Theorem of Thom-Sard (see e.g. \cite{golub}) to the diagram

$$ \begin{array}{lcc}
\Phi_i^{-1} (0) \cap ({\cal U} \times \mC^{(n-i) i})  & \hookrightarrow &
\mC^n \times \mC^{(n-i) i} \\
& \searrow  & \downarrow \\
&  & \mC^{(n-i) i} \end{array} $$

one concludes that the set of all matrices $A^{(i)} \in \mR^{(n-i) i}$ for
which transversality holds is Zariski dense in $\mC^{(n-i) i}$. More
precisely, the affine space $\mQ^{(n-i)i}$ contains a non-empty Zariski open
set of matrices $A^{(i)}$ such that the corresponding coordinate
transformation $(1)$ leads to the desired smoothness of
$\widetilde{W}_i$ in points which are smooth in $W$. \hfill $\Box$ \spar

Let $f \in \mQ[X_1,\ldots,X_n]$ be a non-constant square-free polynomial and
let again $W := \{ x \in \mC^n | f(x) = 0 \}$ be the hypersurface defined
by $f$. Let $\Delta \in \mQ[X_1, \ldots, X_n]$ be the polynomial $\Delta :=
\sum\limits_{j = 1}^{n} \left ( {{\partial f} \over {\partial X_j}}
\right )^2$. Consider the real variety $V := W \cap \mR^n$ and suppose that:

{\em
\begin{itemize}
\item $V$ is non-empty and bounded (and hence compact)
\item the gradient of $f$ is different from zero in all points of $V$\\
(i.e., $V$ is a smooth hypersurface in $\mR^n$ and $f = 0$ is its
regular equation)
\item the variables $X_1, \ldots, X_n$ are in generic position.
\end{itemize}
}

Under these assumption the following problem adapted notion of polar variety
is meaningful and remains consistent with the more general definition of
the same concept (see e.g. \cite{lete}).

\begin{definition}

Let $0\leq i<n$. Consider the linear subspace $X^i$ of  $\mC^n$ 
corresponding to the
linear forms $X_{i+1},\ldots,X_n$, i.e., $X^i := \{ x \in \mC^n |
X_{i+1}(x) = \cdots = X_n(x) = 0 \}$. Then the algebraic subvariety $W_i$ of $\mC^n$ defined as the
Zariski closure of the set 

$$ \{ x \in \mC^n | f(x) = {{\partial f(x)} \over{\partial X_1}} = \cdots
= {{\partial f(x)} \over {\partial X_i}} = 0, \Delta (x) \not= 0 \}$$

is called the (complex) {polar variety\/} of $W$ associated to the
linear subspace $X^i$ of $\mC^n$. The respective real variety is denoted
by $V_i := W_i \cap \mR^n$ and  called the {\rm real polar variety} of $V$
associated to the linear subspace $X^i \cap \mR^n$ of $\mR^n$. Here $W_0$ is
understood to be the Zariski closure of the set $\{ x \in \mC^n ; f(x) = 0,
\Delta(x) \not= 0\}$ and $V_0$ is understood to be $V$.

\end{definition}

\begin{remark} Since by assumption $V$ is a non-empty compact hypersurface
of $\mR^n$ and the variables $X_1, \ldots, X_n$ are in generic position, we
deduce from Proposition $3$ and general considerations on Lagrange
multipliers (as e.g. in \cite{hroy}) or Morse Theory (as e.g. in
\cite{Milnor}) that the real polar variety $V_i$ is non-empty and smooth for
any $0 \le i < n$. In particular, the complex variety $W_i$ is not empty and
the hypersurfaces of $\mC^n$ given by the polynomials $f, {{\partial f}
\over{\partial X_1}}, \ldots, {{\partial f} \over {\partial X_i}}$ intersect
transversally in some dense Zariski open subset of $W_i$ (observe that any
element of $\{ x \in \mC^n ; f(x) = 0, \Delta(x) \not= 0\}$ is a smooth
point of $W$ and apply Proposition 3).  
\end{remark}

Let us observe that the assumption {\em $V$ smooth} implies that the polar variety
$V_i$ can be written as $V_i = \{ x\in \mR^n ; f(x) = {{\partial f(x)}
\over{\partial X_1}} = \cdots {{\partial f(x)} \over {\partial X_i}} = 0
\}$ for any $0\leq i\le n$.

\begin{theorem}

Let $f \in \mQ[X_1,\ldots,X_n]$ be a non-constant square-free polynomial and
let $\Delta := \sum\limits_{j = 1}^{n} \left ( {{\partial f} \over {\partial
X_j}} \right )^2$. Let $W := \{ x \in \mC^n | f(x) = 0 \}$ be the
hypersurface of $\mC^n$ given by the polynomial $f$. Further, suppose that
$V := W \cap \mR^n$ is a non-empty, smooth and bounded hypersurface of
$\mR^n$ with regular equation $f$. Assume that the variables
$X_1, \ldots, X_n$ are in generic position. Finally, let for any $i$, $0
\le i < n$, the complex polar variety $W_i$ of $W$ and the real polar
variety $V_i$ of $V$ be defined as above. With these notations and
assumptions we have~: 

\begin{itemize} 

\item $V_0 \subset W_0 \subset W$, with $W_0 = W$ if and only if $f$ and
$\Delta$ are coprime,

\item $W_i$ is a non-empty equidimensional affine variety of dimension
$n-(i+1)$ being smooth in all its points which are smooth points of $W$,

\item the real part $W_i^\ast $ of the complex polar variety $W_i$ coincides
with the Zariski closure in $\mC^n$ of the real polar variety

$$V_i = \left\{ x \in \mR^n | f(x) = {{\partial f(x)} \over {\partial X_1}}
= \cdots = {{\partial f(x)} \over {\partial X_i}} = 0 \right\} ,$$

\item the ideal $\left(f, {{\partial f} \over {\partial X_1}},
\ldots,{{\partial f} \over {\partial X_i}}\right)_{\Delta}$ is radical.

\end{itemize}
\end{theorem}

{\bf Proof:} 

The first statement is obvious, because $W_0$ is the union of all
irreducible components of $W$ on which $\Delta$ does not vanish identically.
\spar

We show now the second statement. Let $0 \le i < n$ be arbitrarily fixed.
Then the polar variety $W_i$ is non-empty by Remark $5$. Moreover, the
hypersurfaces of $\mC^n$ defined by the polynomials $f, \frac{\partial
f}{\partial X_1}, \ldots, \frac{\partial f}{\partial X_n}$ intersect any
irreducible component of $W_i$ transversally in a non-empty Zariski open
set. This implies that the algebraic variety $W_i$ is a non-empty
equidimensional variety of dimension $n - (i+1)$ and that the polynomials
$f, \frac{\partial f}{\partial X_1}, \ldots, \frac{\partial f}{\partial
X_n}$ form a regular sequence in the ring obtained by localizing $\mQ[X_1,
\ldots, X_n]$ by the polynomials which do not vanish identically on any
irreducible component of $W_i$. More exactly, the polynomials $f,
\frac{\partial f}{\partial X_1}, \ldots, \frac{\partial f}{\partial X_i}$
form a regular sequence in the localized ring $\mQ[X_1, \ldots,
X_n]_{\Delta}$. From Proposition $3$ we deduce that $W_i$ is smooth in all
points which are smooth points of $W$ and that the hypersurfaces of $\mC^n$
defined by the polynomials $f, \frac{\partial f}{\partial X_1}, \ldots,
\frac{\partial f} {\partial X_i}$ intersect transversally in these points.
\spar

Let us show the third statement. The Zariski closure of $V_i$ in $\mC^n$ is
contained in $W_i^\ast$ (this is a simple consequence of the smoothness of
$V_i$). One obtains the reverse inclusion as follows: let $x^\ast \in
W_i^\ast$ be an arbitrary point, and let $C$ be an irreducible component of
$W_i^\ast$ containing this point. Since $C$ is a real component of $W_i$ the
set $C \cap \mR^n $ is not empty and contained in $W_i$. The polar variety
$W_i$ is contained in the algebraic set $\widetilde{W}_i := \left\{ x\in
\mC^n | f(x) = \frac{\partial f (x)}{\partial X_1} = \right.$ $
\left. =\cdots = \frac{\partial
f (x)}{\partial X_i} = 0 \right\}$. Therefore, we have $C \cap V_i \not=
\emptyset$. Moreover, the hypersurfaces of $\mR^n$ defined by the polynomials
$f, \frac{\partial f}{\partial X_1}, \ldots, \frac{\partial f}{\partial
X_i}$ cut out transversally a dense subset of $C \cap V_i$. Thus we
have

\[ \begin{array}{rl} n-(i+1) &= dim_{\R}(C\cap V_i)= dim_{\R} R(C\cap V_i)=\\
&=dim_{\mC} R((C\cap V_i)')\le dim_{\mC} C = n-(i+1). \end{array} \]

(Here $R(C\cap V_i)$ denotes the set of smooth points of $C\cap V_i$ and
$(C\cap V_i)'$ denotes the complexification of $C\cap V_i$.) Thus,
$dim_{\mC} (C \cap V_i)' = dim_{\mC} C = n-(i+1)$ and, hence, $C = (C\cap
V_i)'$. Moreover, $(C\cap V_i)'$ is contained in the Zariski closure of
$V_i$ in $\mC^n$, which implies that $C$ is contained in the Zariski closure
of $V_i$ as well.  \spar

Finally, we show the last statement. Let us consider again the algebraic set

\[ \widetilde{W}_i := \left\{ x\in \mC^n | f(x) = \frac{\partial
f (x)}{\partial X_1} = \cdots = \frac{\partial f (x)}{\partial X_i} = 0
\right\} \]

which contains the polar variety $W_i$. Let $C'$ be any irreducible
component of $W_i$. Then $C'$ is also an irreducible component of
$\widetilde{W}_i$. Moreover, the polynomial $\Delta$ does not vanish
identically on $C'$. By Remark $5$ there exists now a smooth point
$x^\ast$ of $\widetilde{W}_i$ which is contained in $C'$ and in which the
hypersurfaces of $\mC^n$ given by the polynomials $f, \frac{\partial
f}{\partial X_1}, \ldots, \frac{\partial f}{\partial X_i}$ intersect
transversally. \spar

Let $x^\ast = (X_1^\ast, \ldots, X_n^\ast) \in \mC^n$ be fixed in that way. 
Consider the local ring ${\cal O}_{\widetilde{W}_i,x^\ast}$ of the point
$x^\ast$ in the variety $\widetilde{W}_i$ (i.e., ${\cal
O}_{\widetilde{W}_i,x^\ast}$ is the ring of germs of rational functions of
$\widetilde{W}_i$ that are defined in the point $x^\ast$). Algebraically the
local ring ${\cal O}_{\widetilde{W}_i,x^\ast}$ is obtained by dividing the
polynomial ring $\mbox{} \mC[X_1, \ldots, X_n ]$ by the ideal $(f,
\frac{\partial f}{\partial X_1}, \ldots, \frac{\partial f}{\partial X_i})$,
which defines $\widetilde{W}_i$ as an affine variety, and then localizing at
the maximal ideal $(X_1 - X_1^\ast,\ldots ,X_n-X_n^\ast)$ of the point
$x^\ast = (X_1^\ast, \ldots, X_n^\ast)$. Using now standard arguments from
Commutative Algebra and Algebraic Geometry (see e.g. \cite{brod}), one
infers from the fact that the hypersurfaces of $\mC^n$ given by the
polynomials $f, \frac{\partial f}{\partial X_1}, \ldots, \frac{\partial
f}{\partial X_i}$ intersect transversally in $x^\ast$ the conclusion that
${\cal O}_{\widetilde{W}_i,x^\ast}$ is a regular local ring and, hence, an
integral domain.
The fact that ${\cal O}_{\widetilde{W}_i,x^\ast}$ is an integral domain
implies that there exists a uniquely determined irreducible component of
$\widetilde{W}_i$ which contains the smooth point $x^\ast$ (this holds true
for the ordinary, $\mC-$defined Zariski topology as well as for the
$\mQ-$defined one considered here). Therefore, the point $x^\ast$ is
uniquely contained in the irreducible component $C'$ of $\widetilde{W}_i$
(and of $W_i$). 

Since the local ring ${\cal O}_{\widetilde{W}_i,x^\ast}$
is an integral domain, its zero ideal is prime. This implies that the
polynomials $f, \frac{\partial f}{\partial X_1}, \ldots, \frac{\partial
f}{\partial X_i}$ generate a prime ideal in the local ring \linebreak[4] $\mC[X_1, \ldots
,X_n]_{(X_1 - X_1^\ast,\ldots ,X_n-X_n^\ast)}$. Hence, the isolated
primary component of the polynomial ideal \linebreak[4] $(f, \frac{\partial f}{\partial
X_1}, \ldots, \frac{\partial f}{\partial X_i})$ in $\mQ[X_1, \ldots, X_n]$,
which corresponds to the irreducible component $C'$, is itself a prime ideal.
Since this is true for any irreducible component of $W_i$ and since $W_i$
defined by discarding from $\widetilde{W}_i$ the irreducible components
contained in the hypersurface of $\mC^n$ given by the polynomial
$\Delta$, we conclude that the ideal $(f, \frac{\partial f}{\partial X_1},
\ldots, \frac{\partial f}{\partial X_i})_{\Delta}$ of $\mQ[X_1, \ldots,
X_n]_{\Delta}$ is an intersection of prime ideals and, hence, radical. This
completes the proof of Theorem 6.\hfill $\Box$

\begin{remark}\label{rem7}
Under the assumptions of Theorem 6, we observe that for any $i$,  $0 \le i < n$,
the following inclusions hold among the different non-empty varieties
introduced up to now, namely

\[ V_i \subset V\mbox { and } V_i \subset W^\ast_i \subset W_i \subset
\widetilde{W}_i.\]

Here $V$ is the bounded and smooth real hypersurface we consider in this
paper, $W_i$ and $V_i$ are the polar varieties introduced in Definition 4,
$W_i^\ast$ is the real part of $W_i$ according to Definition 1, and
$\widetilde{W}_i$ is the complex affine variety introduced in the proof of
Theorem 6. With respect to Theorem 6 our settings and assumptions imply that
$n-(i+1) = dim_{\mC} W_i = dim_{\mC} W^\ast_i = dim_{\R} V_i $ holds.  By
our smoothness assumption and the generic choice of the variables we have
for the respective sets of smooth points the inclusions:

\[ V_i = R(V_i)\subset R(W_i) \subset R(\widetilde{W}_i) \subset R(W) \]

(Here $W$ is the affine hypersurface $W = \{x \in \mC^n | f(x) = 0\}$ of
$\mC^n$.)
 
\end{remark}
\end{section}

\begin{section}{Algorithms and Complexity}

The preceeding study of adapted polar varieties enables us to state our 
first complexity result:

\begin{theorem}\label{theo8}

Let $n, d, \delta, L$ be natural numbers. Then there exists an arithmetic
network $\cal{N}$ over $\mQ$ of size $(nd\delta L)^{O(1)}$ with the
following properties:
\spar

Let  $f \in \mQ[X_1, \ldots, X_n]$ be a non-constant polynomial of degree at
most $d$ and suppose that $f$ is given by a division--free straight--line
program $\beta$ in $\mQ[X_1, \ldots, X_n]$ of length at most $L$. Let
$\Delta := \sum\limits_{j = 1}^{n} \left ( {{\partial f} \over {\partial
X_j}} \right )^2$, $W := \{x \in \mC^n | f(x) = 0\}$, $V := W \cap \mR^n =
\{x \in \mR^n | f(x) = 0\}$ and suppose that the variables $X_1, \ldots,
X_n$ are in "sufficiently generic" position. For $0 \le i <n$ let $W_i$ be
the Zariski closure in $\mC^n$ of the set

$$\{x \in \mC^n | f(x) =  \frac{\partial f (x) }{\partial X_1} = \cdots = 
\frac{\partial f (x)}{\partial X_i} = 0, \Delta(x) \not= 0\}$$

(thus $W_i$ is the polar variety of $W$ associated to the linear space $X^i
= \{ x \in \mC^n | X_{i+1}(x) = \cdots = X_n(x) = 0 \}$ according to
Definition 4). Let $\delta_i := deg W_i$ be  the geometric degree of $W_i$ and
assume that $\delta \ge max \{ \delta_i | 1 \le i <n\}$ holds. 
\spar

The algorithm represented by the arithmetic network $\cal{N}$ starts from
the straight--line program $\beta$ as input and decides first whether the
complex algebraic variety $W_{n-1}$ is zero--dimensional. If this is the
case the network $\cal{N}$ produces a straight--line program of length $(nd
\delta L)^{O(1)}$ in $\mQ$ which represents the coefficients of $n+1$
univariate polynomials $q,p_1, \ldots, p_n \in \mQ[X_n]$ satisfying the
following conditions:

\begin{enumerate}
\item $deg (q) = \delta_{n-1} = deg W_{n-1}$
\item $max \{ deg (p_i) | 1 \le i \le n \} < \delta_{n-1}$
\item $W_{n-1} = \{ (p_1(u), \ldots, p_n(u)) | u \in \mC, q(u) = 0\}$.
\end{enumerate}

Moreover, the algorithm represented by the arithmetic network $\cal{N}$
decides whether the semialgebraic set $W_{n-1} \cap \mR^n$ is non-empty. If
this is the case the network $\cal{N}$ produces not more than $\delta_{n-1}$
sign sequences of $\{ -1, 0, 1\}^{\delta_{n-1}}$ which codify the real
zeros of $q$ "\`a la Thom" (\cite{CosteRoy}). In this way, $\cal{N}$
describes the non-empty finite set $W_{n-1} \cap \mR^n$.

\end{theorem}

>From the output of this algorithm we may deduce the following information:

\begin{itemize}

\item If the complex variety $W_{n-1}$ is not zero-dimensional or if
$W_{n-1}$ is zero-dimensional and $W_{n-1} \cap \mR^n$ is empty we conclude
that $V$ is not a compact smooth hypersurface of $\mR^n$ with regular
equation $f$.

\item  If $V$ is a compact smooth hypersurface of $\mR^n$ with regular
equation $f$, then $W_{n-1} \cap \mR^n$ is non-empty and contains for any
connected component of $V$ at least one point which the network $\cal{N}$
codifies \`a la Thom as a real zero of the polynomial $q$.

\end{itemize}

\begin{remark}\label{rem9}

The hypothesis that the variables $X_1, \ldots ,X_n$ are in "sufficiently
generic" position is not really restrictive since any $\mQ-$linear
coordinate change increases the length of the input straight-line program
$\beta$ only by an unessential additive term of $O(n^3)$. Moreover, by
\cite{HeSchn} Theorem 4.4, any genericity condition which the algorithm
might require can be satisfied by adding to the arithmetic network ${\cal N}$ an
extra number of nodes which is polynomial in the input parameters $n, d,
\delta, L$.  \end{remark}

\begin{remark}\label{rem10} From the B\'ezout Theorem we deduce the estimation
$max\{ \delta_i | 0\leq i < n\} \leq d(d-1)^{n-1}<d^n$. Moreover $f$ can
always be evaluated by a division-free straight--line program in
$\mQ[X_1,\ldots,X_n]$ of length $d^n$. Thus fixing $\delta:=d(d-1)^{n-1}$ and
$L:=d^n$ one is concerned with a worst case situation in which the
statement of Theorem 8 just reproduces the main complexity results of
\cite{grivo1}, \cite{Griegorev87}, \cite{hroy}, \cite{hroy1}, \cite{hroy2},
\cite{canny}, \cite{Renegar1}, \cite{rene}, \cite{basu} in case of a compact smooth
hypersurface of $\mR^n$ given by a regular equation of degree $d$. The
interest in Theorem 8 lies in the fact that $\delta$ may be much smaller
than the "B\'ezout number" $d(d-1)^{n-1}$ and $L$ smaller than $d^n$ in many
concrete and interesting cases.

\end{remark}

{\bf Proof of Theorem 8:} 

Since by \cite{bast} and \cite{morg} we may derive the straight-line program
$\beta$ representing the polynomial $f$ in time linear in $L$, we may
suppose without loss of generality that $\beta$ represents also the
polynomial $\Delta$. Applying now the algorithm underlying \cite{gh2}
Proposition 18 together with the modifications introduced by \cite{gh3}
Theorem $28$ (compare also \cite{gh3} Theorem 16 and its proof), we find an
arithmetic Network $\cal{N'}$ with parameters in $\mQ$ of size $(nd \delta
L)^{O(1)}$ which decides whether the polynomials $f, \frac{\partial
f}{\partial X_1},\ldots, \frac{\partial f}{\partial X_{n-1}}$ form a secant
family avoiding the hypersurface of $\mC^n$ defined by the polynomial
$\Delta$. This is exactly the case if $W_{n-1}$ ist zero-dimensional. \spar

Suppose now that the polynomials $(f, \frac{\partial f}{\partial
X_1},
\ldots, \frac{\partial f}{\partial X_{n-1}})$ form such a secant family.
Then the arithmetic network $\cal{N'}$ which we obtained before applying
\cite{gh2} Proposition $18$ and \cite{gh3} Theorem $28$ to the input
$f, \frac{\partial f}{\partial X_1}, \ldots, \frac{\partial f}{\partial
X_{n-1}}$ and $\Delta$ produces a straight-line program in $\mQ$ which
represents the coefficients of polynomials $q, p_1, \ldots, p_n \in \mQ[X_n]$ 
characterizing the part $W_{n-1}$ of the complex variety 
$\widetilde{W}_{n-1} := \left\{ x\in \mC^n | f(x) = \frac{\partial f(x)}{\partial
X_1} = \cdots = \frac{\partial f(x)}{\partial X_{n-1}} = 0 \right\}$
which avoids the hypersurface $\{ x \in \mC^n | \Delta (x) = 0 \}$.  More
precisely, the output $q, p_1, \ldots, p_n$ of the network ${\cal N'}$ satisfies
the conditions $(1), (2), (3)$ in the statement of Theorem $8$.  \spar

Now applying for example the main (i.e., the only correct) algorithm of
\cite{BenOrKozenReif} (see also \cite{RoySzpirglas} for refinements) by
adding suitable comparison gates for positiveness of rational numbers, we
may extend $\cal{N'}$ to an arithmetic network $\cal{N}$ of asymptotically
the same size $(nd \delta L)^{O(1)}$, which decides whether the polynomial
$q$ has any real zero. Moreover, without loss of generality the arithmetic
network $\cal{N}$ codifies any existing zeros of $q$ \`a la Thom (see
\cite{CosteRoy}, \cite{RoySzpirglas}). From general considerations of Morse
Theory (see e.g. \cite{Milnor}) or more elementary from the results and
techniques of \cite{hroy}, \cite{hroy2} one sees that in the case where $f$
is a regular equation of a bounded smooth hypersurface $V$ of $\mR^n$, the
arithmetic network $\cal{N}$ codifies for each connected component of $V$ at
least one representative point. This finishes the proof of Theorem $8$. 
\hfill $\Box$ \spar

Roughly speaking the arithmetic network $\cal{N}$ of Theorem $8$ decides
whether a given polynomial $f \in \mQ[X_1, \ldots, X_n]$ is a regular
equation of a bounded (i.e., compact) smooth hypersurface $V$ of $\mR^n$. If
this is the case $\cal{N}$ computes for any connected component of $V$ at
least one representative point. The size of $\cal{N}$ depends polynomially
on the number of variables $n$, the degree $d$ and the straight-line program
complexity $L$ of $f$ and finally on the degree $\delta$ of certain complex
polar varieties $W_i$ associated to the equation $f$. 

The nature of the answer the network $\cal{N}$ gives us about the
algorithmic problem is satisfactory. However, this is not the case for the
size of $\cal{N}$, which measures the complexity of the underlying
algorithm, since this complexity depends on the parameter $\delta$
being related rather to the complex considerations  than to the real 
ones. We are going now to describe a procedure whose complexity is
polynomial only in the {\em real\/} degree of the polar varieties $W_i$
instead of their complex degree. The theoretical (not necessarily the
practical) price we have to pay for this complexity improvement is
relatively high:

\begin{itemize}

\item our new procedure does not decide any more whether the input polynomial
is a regular equation of a bounded smooth hypersurface $V$ of $\mR^n$. We
have to assume that this is already known. Therefore the new algorithm can
only be used in order to {\em solve\/} the real equation $f = 0$, but not to
decide its consistency ({\em solving\/} means here that the algorithm
produces at least one representative point for each connected component of
$V$).

\item our new algorithm requires the support of the following two external
subroutines whose theoretical complexity estimates are not really taken into
account here although their practical complexity may be considered as
"polynomial":

\begin{itemize}

\item the first subroutine we need is a factorization algorithm for
univariate polynomials over $\mQ$. In the bit complexity model the problem
of factorizing univariate polynomials over $\mQ$ is known to be
polynomial (\cite{LLL}), whereas in the arithmetic model we
are considering here this question is more intricate (\cite{vzGaSe91}).
In the extended complexity model we are going to consider, the cost of
factorizing a univariate polynomial of degree $D$ over $\mQ$, (given by its
coefficients) is accounted as $D^{O(1)}$.

\item  the second subroutine allows us to discard non-real irreducible
components of the occuring complex polar varieties. This second subroutine
starts from a straight-line program for a single polynomial in $\mQ[X_1,
\ldots, X_n]$ as input and decides whether this polynomial has a real zero
(however without actually exhibiting it if there is one). Again this
subroutine is taken into account at polynomial cost.

\end{itemize}
 
\item We call an arithmetic network  over $\mQ$ {\em extended\/} if it
contains extra nodes corresponding to the first and second subroutine.

\end{itemize}

Fix for the moment natural numbers $n, d, \delta^\ast$ and $L$. We suppose
that a division-free straight-line program $\beta$ in $\mQ[X_1, \ldots,
X_n]$ of length at most $L$ is given such that $\beta$ represents a
non-constant polynomial $f \in \mQ[X_1, \ldots, X_n]$ of degree at most $d$.
Let again $\Delta := \sum \limits_{j=1}^n \left ( \frac{\partial f}{\partial
X_j} \right ) ^2$ and suppose that $f$ is a regular equation of a
(non-empty) bounded smooth hypersurface $V$ of $\mR^n$. Let $W:=\{ x\in\mC^n
; f(x) =0\}$ be the complex hypersurface of $\mC^n$ defined by the
polynomial $f$ and suppose that the variables $X_1,\ldots,X_n$ are in
generic position. Fix $0 \le i <n$ arbitrarily. Let as in Definition $4$ the
complex variety $W_i$ be the the Zariski closure in $\mC^n$ of the set

$$ \{ x \in \mC^n | f(x) = {{\partial f(x)} \over{\partial X_1}} = \cdots
= {{\partial f(x)} \over {\partial X_i}} = 0, \Delta (x) \not= 0 \} ,$$

i.e., $W_i$ is the polar variety of the complex hypersurface $W$ associated
to the linear subspace $X^i := \{x \in \mC^n | X_{i+1}(x) =0, \ldots, X_n(x)
= 0\}$.
\spar

Let $V_i := W_i \cap \mR^n$ the corresponding polar variety of the real
hypersurface $V$. Let $\delta_i^\ast$ be the real degree of the polar
variety $W_i,$ i.e., the geometric degree of $W_i^\ast$ (see Definition $1$).
By Theorem $6$ the quantity $\delta_i^\ast$ is also the geometric degree of
the Zariski closure in $\mC^n$ of the real variety $V_i$, i.e., of the
complexification of $V_i$. Let $r:=n-(i+1)$. Since the variables $X_1,
\ldots, X_n$ are in generic position with respect to all our geometric data,
they are also in Noether position with respect to the complex variety $W_i$,
the variables $X_1, \ldots, X_r$ being free (see \cite{gh1}, \cite{gh2} for
details). Finally, suppose that $\delta^\ast \ge max \{\delta_i^\ast |  0
\le i <n \}$ holds. \spar

With these notations and assumptions, we have the following {\em real\/}
version of \cite{gh2} Proposition $17$:

\begin{lemma}\label{lem11} 

Let $n, d, \delta^\ast L$ be given natural numbers as before and fix $0 \le
i <n$ and $r:=n-(i+1)$. Then there exists an extended arithmetic network
$\cal{N}$ with parameters in $\mQ$ of size $(i d \delta^\ast L)^{O(1)}$
which for any non-constant polynomial $f \in \mQ[X_1, \ldots, X_n]$
satisfying the assumptions above produces a division-free straight-line
program $\beta_i$ in $\mQ[X_1, \ldots, X_r]$ such that $\beta_i$ represents
a non-zero polynomial $\varrho \in \mQ[X_1, \ldots, X_r]$ and the
coefficients with respect to $X_{r+1}$ of certain polynomials $q, p_1,
\ldots, p_n \in \mQ[X_1, \ldots, X_{r+1}]$ having the following properties:

\begin{itemize}

\item[(i)] the polynomial $q$ is monic and separable  in $X_{r+1}$, square-free
and its degree satisfies $deg q = deg_{X_{r+1}} q = \delta^\ast_{i}= deg
W^\ast_{i} \le \delta^\ast$,

\item[(ii)] the polynomial $\varrho$ is the discriminant of $q$ with respect
to the variable $X_{r+1}$ and its degree can be estimated as $deg \varrho
\le 2(\delta_i^\ast)^3$,

\item[(iii)] the polynomials $p_1, \ldots, p_n$ satisfy the degree
bounds

$max \{ deg_{X_{r+1}} p_k | 1\leq k \le n \} < \delta^\ast_{i}$,\ \ \ 
$max \{ deg p_k | 1 \leq k \le n \} = 2( \delta_i^\ast)^3,$

\item[(iv)] the ideal $(q, \varrho X_1-p_1, \ldots ,\varrho
X_n-p_n)_\varrho$ generated by the polynomials $q, \varrho X_1-p_1, \ldots,
\varrho X_n-p_n$ in the localization $\mQ[X_1,\ldots ,X_n]_\varrho$ is the
vanishing ideal of the affine variety $(W_i^\ast)_\varrho := \{ x \in
W_i^\ast | \varrho (x) \not= 0 \}$. Moreover, $(W_i^\ast)_\varrho$ is a dense
Zariski open subset of the complex variety $W_i^\ast$.

\item[(v)] the length of the straight-line program $\beta_i$ is of the order
$(id\delta^\ast L)^{O(1)}$.

\end{itemize}

\end{lemma}

{\bf Proof:} 

The proof of this lemma follows the general lines of the proof of Theorem
$8$ and is again based on the algorithm underlying \cite{gh2} Proposition
$17$ together with the modifications introduced by \cite{gh3} Theorem $28$.
\spar

First we observe that by \cite{bast} and \cite{morg} we are able to derive
the straight-line program $\beta$ representing the input polynomial $f$ at
cost linear in $L$. Thus we may suppose without loss of generality that
$\beta$ represents both $f$ and $\Delta$. \spar

We show Lemma 11 by the exhibition of a recursive procedure in $0 \le i <n$
under the assumption that the first and second subroutine as introduced
before are available. First put $i := 0$ and let $\beta_0$ be the
straight-line program $\beta$ which represents $f$ and $\Delta$. Since the
variables $X_1, \ldots, X_n$ are in generic position, the polynomials $f$
and $\Delta$ are monic with respect to the variable $X_n$ and satisfy the
conditions $d \ge deg f = deg_{X_n} f$ and $2d \ge deg \Delta =
deg_{X_n} \Delta$. \spar

Let $R_0 := \mQ[X_1, \ldots, X_{n-1}]$ and consider $f$ and $\Delta$ as
univariate polynomials in $X_n$ with coefficients in $R_0$. Recall that they
are monic. Interpolating them in $2d+1$ arbitrarily chosen distinct rational
points, we obtain a division-free straight-line program in $R_0 = \mQ[X_1,
\ldots, X_n]$ which represents the coefficients of $f$ and $\Delta$ with
respect to $X_n$. This straight-line program has length $L d^{O(1)}$.
\spar

We apply now \cite{gh2} Lemma $8$ in order to obtain the maximal common
divisor of $f$ and $\Delta$ which is again a monic polynomial in $R_0[X_n]$
which we may suppose to be given by a division-free straight-line program in
$R_0 = \mQ[X_1, \ldots, X_{n-1}]$ representing its coefficients with respect
to $X_n$. Dividing $f$ by this maximal common divisor in $R_0[X_n]$ as in
the Noether Normalization procedure in \cite{gh2}, we obtain a polynomial
$\bar{q} \in R_0 [X_n] = \mQ[X_1, \ldots, X_n]$ whose coefficients with
respect to the variable $X_n$ are represented by a division-free
straight-line program $\bar{\beta}_1$ in $\mQ[X_1, \ldots, X_{n-1}]$. The
polynomial $\bar{q}$ is monic in $X_n$, it is square-free and it is a divisor
of $f$. Moreover we have $W_0 = \{ x \in \mC^n | \bar{q} (x) = 0 \}$
and $\bar{q}$ is the minimal polynomial of the hypersurface $W_0$ of
$\mC^n$. The degree of the polynomial $\bar q$ satisfies the condition $deg
\bar{q} = deg_{X_0} \bar{q} = deg W_0 $. \spar

The straight-line program $\bar{\beta}_1$ which represents the coefficients
of $\bar{q}$ with respect to the variable $X_n$ has length $(dL)^{O(1)}$. In
order to finish the recursive construction for the case $i:=0$ it is
sufficient to find the factor $q$ of $\bar{q}$ which defines the real part
$W_0^\ast$ of $W_0$. For this purpose we consider the projection map $\mC^n
\rightarrow \mC^{n-1}$ which maps each point of $\mC^n$ onto its first $n-1$
coordinates. Since the variables $X_1, \ldots, X_n$ are in generic position,
the projection map induces a finite surjective morphism $\pi : W_0
\rightarrow \mC^{n-1}$. We choose a generic {\em lifting point\/} $t = (t_1,
\ldots, t_{n-1}) \in \mQ^{n-1}$ with rational coordinates $t_1, \ldots, t_{n-1}$
(see \cite{gh3} for the notion of {\em lifting point}). Observe that the
irreducible components of $W_0$ are the hypersurfaces of $\mC^n$ defined by
the $\mQ$-irreducible factors of $\bar{q}$ which we denote by
$q_1, \ldots, q_s$.  \spar

Without loss of generality we may assume that for $1 \le m \le s$ the
irreducible polynomials $q_1, \ldots, q_m$ define the real irreducible
components of $W_0$. Thus it is clear that the factor $q$ of $\bar{q}$
we are looking for is $q := q_1\cdots q_m$. It suffices therefore
to find all irreducible factors $q_1, \ldots, q_s$ of $\bar{q}$ and
then to discard the factors $q_{m+1}, \ldots, q_s$. \spar

In order to find the polynomials $q_1, \ldots, q_s$, we specialize the
variables $X_1, \ldots,$ $\ldots X_{n-1}$ into the coordinates $t_1$,$ \ldots $, $t_{n-1}$
of the rational point $t\in\mQ^{n-1}$.  We obtain thus the univariate
polynomial $\bar{q} (t, X_n) := \bar{q} (t_1, \ldots, t_{n-1}, X_n)$ $\in
\mQ[X_n]$ which decomposes into $\bar{q} (t, X_n) = {q}_1 (t, X_n) \cdots
 {q}_s (t, X_n)$ in $\mQ[X_n]$. Since the lifting point $t$ was
chosen generically in $\mQ^{n-1}$, Hilbert's Irreducibility Theorem (see
\cite{lang}) implies that the polynomials ${q}_1 (t, X_n), \ldots, {q}_s (t,
X_n)$ are irreducible over $\mQ$. Specializing the variables
$X_1,\ldots,X_{n-1}$ in the straight-line program $\bar{\beta_1}$ into the
values $t_1, \ldots, t_{n-1}$ we obtain an arithmetic circuit in $\mQ$ which
represents the coefficients of $\bar{q}(t,X_n)$. By a call to the first
subroutine we obtain the coefficients of the polynomials ${q}_1(t,
X_n), \ldots, {q}_s(t, X_n)$. Applying to these polynomials the
lifting procedure which we are going to explain below in a slightly more
general context, we find a division-free straight-line program in $\mQ[X_1,
\ldots, X_{n-1}]$ of size $(dL)^{O(1)}$ which represents the coefficients of
the polynomials $q_1, \ldots, q_s$ with respect to the variable $X_n$.\\
In order to finish the case  $i=0$ we have to
identify algorithmically the polynomials $q_1,\ldots ,q_m$ that
define the irreducible real components of $W_0$ and, hence, those
of $W^\ast$. Then, the product $q=q_1 \cdots q_m$ is easily
obtained.  Observe that $q$ is the minimal polynomial of the
hypersurface $W_0$. From the assumption that $V=W\cap \R^n$ is
a smooth real hypersurface one deduces that $V_0 = W_0^\ast \cap
\R^n = W_0 \cap \R^n$ holds. Since $f$ is a regular equation of
$V$ and since the polynomials $\bar{q}$ and $q$ are factors of
$f$, one sees immediately that $\bar{q}$ and $q$ are also
regular equations of $V$.  This implies that each of the
polynomials $q_1,\ldots ,q_s$ admitting a real zero $x \in \R$ has a
non-vanishing gradient in $x$. Thus, any polynomial of $q_1, \ldots ,
q_s$ admitting a real zero belongs to $q_1,\ldots, q_m$. Hence,
by a call to the second subroutine, we are able to find the
polynomials $q_1,\ldots ,q_m$, and therefore the polynomial
$q=q_1\cdots q_m$.\\ 
Now we extend the division-free
straight-line program representing the polynomials
$q_1,\ldots , q_s$ to a circuit in $\mQ [X_1,\ldots ,X_n]$ of
size $(dL)^{O(1)}$ which computes the polynomial $q=q_1\cdots
q_s$. Interpolating $q$ in the variable $X_n$ as before, this
circuit provides a division-free straight-line program $\beta_1$
in $\mQ [ X_1, \ldots , X_n]$ of size $(dL)^{O(1)}$ which
represents the coefficients of $q$ with respect to the variable
$X_n$. Without changing its order of complexity we extend
$\beta_1$ to a division-free circuit in $\mQ [X_1,\ldots
,X_{n-1}]$ that computes also the discriminant $\varrho$ of $q$
with respect to the variable $X_n$ and the polynomials $\varrho
X_1,\ldots ,\varrho X_{n-1}$.\\ 
Let $p_1 := \varrho X_1,\ldots ,p_{n-1}:= \varrho
X_{n-1}, p_n := \varrho X_n \in \mQ [X_1, \ldots , X_{n-1}, X_n]$.
One sees immediately that the polynomials $\varrho \in \mQ [
X_1,\ldots , X_{n-1}]$ and $q, p_1, \ldots , p_n \in \mQ [
X_1,\ldots , X_{n-1}, X_n]$ satisfy the conditions (i) - (iv) of
Lemma 11 for $i=0$.  Furthermore, $\beta_1$ is a division-free
straight-line program in $\mQ [ X_1, \ldots ,X_{n-1}]$ of size
$(dL)^{O(1)}$ which computes $\varrho$ and the coefficients of
$q, p_1, \ldots , p_n$ with respect to the variable $X_n$.
By construction the output circuit $\beta_1$ can be produced
from the input circuit $\beta$ by an extended arithmetic network
over $\mQ$ of size $(dL)^{O(1)}$. This finishes the description
of the first stage in our recursive procedure.\\

We consider now the case of $0 < i \le n$ and set $ r := n - (i+1)$.
Suppose that there is given a division-free straight-line program $\beta _{i-1}$
in $\mQ [ X_1,\ldots , X_{r+1}]$ of size $\Lambda_{i-1}$ that reepresents
a non-zero polynomial $ \varrho ' \in \mQ [ X_1,\ldots ,X_{r+1}]$ and the 
coefficients with respect to $X_{r+2}$ of certain polynomials $q', p'_1,
\ldots , p'_n \in \mQ [X_1,\ldots , X_{r+1},X_{r+2}]$.
These polynomials have the following properties: $q'$ is monic and separable in
$X_{r+2}$ and satisfies the degree condition $deg q' = deg_{X_{r+1}} q' = 
\delta_{i-1}^\ast $, $\varrho '$ is the discriminant of $q'$ 
with respect ot $X_{r+2}$,
the polynomials $p'_1,\ldots , p'_n$ satisfy the degree bound
$max\{ deg_{X_{r+2}} p'_k | 1 \le k \le n\} < \delta^\ast_{i-1}$ and the ideal
$(q', \varrho ' X_1-p'_1,\ldots , \varrho ' X_n-p'_n)_{\varrho '}$ of the
localized ring $\mQ [ X_1,\ldots , X_n]_\varrho $ is the vanishing ideal of the 
affine variety $(W_{i-1}^\ast)_{\varrho '}$. Observe that  $(W_{i-1}^\ast)_{\varrho '}$
is a Zariski open dense subset of  $(W_{i-1}^\ast)$. Let $Z$ be the Zariski 
closure in $\mC^n$ of $\{ x \in W^\ast_{i-1} | \frac{\partial f (x)}{\partial  X_i}
= 0, \Delta (x) \not= 0 \}$. We have $W^\ast_i \subset Z \subset W_i$ and 
$Z$ is at least the union of all real irreducible components of $W_i$.
In particular, all irreducible components of $Z$ are irreducible components of 
$W_i$. Moreover, we have $deg  Z \le d \delta^\ast_{i-1}$.\\

Now we apply the procedure underlying \cite{gh2}, Proposition 15, to the 
straight-line programs $\beta_{i-1}$ and $\beta$ representing the polynomials
$\varrho ', q', p'_1, \ldots , p'_n$ and $\frac{\partial f}{\partial X_i}$
in order to produce an explicit description of the algebraic variety
$\{ x \in W^\ast_{i-1} | \frac{\partial f (x)}{\partial X_i} = 0 \}$. By means
of the algorithm of \cite{gh2}, Subsection 5.1.3, we clear out the irreducible 
components of this variety contained in the hypersurface $\{ x \in \mC^n | 
\Delta (x) = 0\}$.
Finally, we obtain a division-free straight-line program $\bar{\mu}$ in $\mQ
[ X_1, \ldots , X_r ]$ of size $i(L+\Lambda_{i-1})(d \delta^\ast_{i-1})$ which
represents a non-zero polynomial $ \bar{\varrho} \in \mQ [X_1,\ldots ,X_r]$ and the 
coefficients with respect to $X_{r+1}$ of certain polynomials  $\bar{q}, \bar{p}_1,
\ldots , \bar{p}_n \in \mQ [X_1, \ldots , X_{r+1}]$. 
The latter polynomials have the following properties: $\bar{q}$ is monic and 
separable with respect to  $X_{r+1}$ and satisfies the degree condition 
$deg \bar{q} = deg_{X_{r+1}} \bar{q} = deg Z$, $\bar{\varrho}$ is the discriminant
of $\bar{q}$ with respect to $X_{r+1}$, the polynomials $\bar{p}_1, \ldots , \bar{p}_n$
satisfy the degree bound $max \{ deg_{X_{r+1}} \bar{p}_k | 1 \le k \le n \} <
deg Z$, and the ideal $ (\bar{q}, \bar{\varrho} X_1 - \bar{p}_1, \ldots , 
\bar{\varrho}X_n - \bar{p}_n)_{\bar{\varrho}}$ of the localized ring $\mQ [ X_1,
\ldots , X_n]_{\bar{\varrho}}$ is the vanishing ideal of the affine variety
$Z_{\bar{\varrho}}$. Observe again that $Z_{\bar{\varrho}}$ is a Zariski
open dense subset of $Z$. By \cite{gh2}, Proposition 15 and Subsection 5.1.3,
there exists an arithmetic network ${\cal N}_i$ with parameters in $\mQ$ 
which produces from the input circuits $\beta_{i-1}$ and $\beta$ the output
circuit $\bar{\mu}$ and has size $i(d \delta^\ast_{i-1} L \Lambda_{i-1})^{O(1)}$.\\
Let $q_1,\ldots ,q_s \in \mQ [ X_1,\ldots ,X_{r+1}]$ be the $\mQ$-irreducible 
factors of $\bar{q}$. Since $\bar{q}$ is monic and separable in $X_{r+1}$, we have
$ \bar{q} = q_1\cdots  q_s$. From the assumption that the variables $X_1,
\ldots , X_n$ are in generic position we deduce that each irreducible component
of the algebraic variety $Z$ is represented by exactly one of the irreducible
polynomials $q_1,\ldots ,q_s$. This means that $Z$ has $s$ irreducible components,
say $C_1,\ldots ,C_s$, such that for $1\le l\le s$ the irreducible component $C_l$
is identical with the Zariski closure in $\mC^n$ of the set
\bean
\{ x = (X_1,\ldots , X_n)\in \mC^n  | &\bar{\varrho}(X_1,\ldots , X_r)X_1-
\bar{p}_1(X_1,\ldots , X_{r+1})= 0,\\
&\ldots \ldots \ldots \ldots \ldots \ldots \ldots 
\ldots \ldots \ldots \ldots \ldots \\
&\bar{\varrho}(X_1,\ldots , X_r)X_n-\bar{p}_n(X_1,\ldots , X_{r+1})= 0,\\
&\!\!\!q_l (X_1,\ldots , X_r)= 0, \bar{\varrho}(X_1,\ldots , X_r)\not= 0\}.
\eean
Now suppose that the real irreducible components of $Z$ (and hence, the one of 
$W_i$) are represented in this way by the polynomials $q_1,\ldots ,q_m$ and let
$q:= q_1 \cdots q_m$. The polynomial $q$ is monic and separable in $X_{r+1}$
and satisfies the degree condition $deg q = deg_{X_{r+1}} q = \delta^\ast_i$. 
Moreover, we have $W^\ast_i =C_1 \cup \ldots \cup C_m$.\\

Now we try to find a straight-line program $\mu$ whose length is of order
$(i d  \delta^\ast_i  L)^{O(1)}$ (hence, independent of the length 
$\Lambda_{i-1}$ of the circuit $\beta_{i-1}$) and which represents the coefficients
of the polynomials $q_1,\ldots , q_s$ and, finally, the polynomial $q$. Adding
to the arithmetic network ${\cal N}_i$ order of $(i d \delta^\ast_i  L
\Lambda_{i-1})^{O(1)}$ extra nodes we find as in the proof of \cite{gh3}, 
Proposition  30, a " sufficiently generic" lifting point $t = (t_1,\ldots ,
t_r)\in \mQ^r$ for the algebraic variety $Z$ (see \cite{gh3}, Definition 19,
for the notion of a lifting point). By the generic choice of the point $t$ we 
deduce from Hilbert's Irreducibility Theorem that $ q_1(t, X_{r+1}), \ldots , q_s
(t, X_{r+1})$ are irreducible polynomials of $\mQ [X_{r+1}]$. Thus the identity
$\bar{q}(t,X_{r+1})=q_1(t,X_{r+1}) \cdots q_s(t,X_{r+1})$ represents the 
decomposition of the polynomial $\bar{q}(t,X_{r+1})\in \mQ [X_{r+1}]$ into its
irreducible factors.\\
Specializing in the straight-line program $\bar{\mu}$ the variables $X_1, \ldots ,
X_r$ into the values $t_1,\ldots ,t_r$ we obtain an arithmetic circuit in $\mQ$
that represents the coefficients of the univariate polynomial $\bar{q}(t,X_{r+1})$.
Adding to the arithmetric network, without changing its asymptotical size, some
extra nodes we may suppose that ${\cal N}_i$ represents the non-zero rational number
$\bar{\varrho}(t)$ and the coefficients of the univariate polynomials
$\bar{q}(t,X_{r+1}), \bar{p}_1(t,X_{r+1}), \ldots , \bar{p}_n(t,X_{r+1})$.
Observe that $deg \bar{q} (t,X_{r+1}) = deg _{X_{r+1}}\bar{q} = deg \bar{q} \le
d \delta^\ast_{i-1}$ holds. Therefore we are able to find the irreducible 
factors $q_1(t,X_{r+1}),\ldots ,q_s(t,X_{r+1})$ of $ \bar{q}(t,X_{r+1})$ by a call
to the first subroutine at a supplementary cost of $(d \delta^\ast_{i-1})^{O(1)}$.
Adding to the arithmetic network ${\cal N}_i$, without changing its asymptotical 
complexity, some extra nodes we may suppose that ${\cal N}_i$ represents for each
$1\le l\le s$ the rational number $\bar{\varrho}(t)$ and the coefficients of the
polynomials $q_l(t,X_{r+1}), \bar{p}_1(t,X_{r+1}), \ldots , \bar{p}_n
(t,X_{r+1}).$ Observe that ${\cal N}_i$ is now an {\em extended} arithmetic network.
For a fixed $l$, $1\le l \le s,$ the set $C_l \cap (\{ t\} \times \mC^{n-r})$ 
is  the lifting fiber of the point $t$ in the irreducible component $C_l$ of
$Z$.The polynomials $q_l(t,X_{r+1}), \frac{1}{\bar{\varrho} (t)} \bar{p}_1
(t,X_{r+1}), \ldots , \frac{1}{\bar{\varrho}(t)} \bar{p}_n (t,X_{r+1})$ 
represent a geometric solution of this lifting fiber.
This means that the identity
\[
C_l\cap (\{ t\} \times  \mC^{n-r})= \left\{ \left( 
\frac{\bar{p}_1(t,u)}{\bar{\varrho}(t)} , \ldots ,
\frac{\bar{p}_n(t,u)}{\bar{\varrho}(t)} \right) | u \in \mC, \; q_l (t,u) =
0 \right\} 
\]
holds.\\

Applying the algorithm underlying \cite{gh3}, Theorem 28, to the input $\beta$,
$t = (t_1,\ldots t_r)$, $\bar{\varrho} (t),$ $q_l(t, X_{r+1}),
 \bar{p}_1(t,X_{r+1}),
\ldots , \bar{p}_n (t,X_{r+1})$ we obtain a division-free straight-line program
in \linebreak[4] $\mQ [X_1, \ldots , X_r]$ having a length of order $(i d deg C_l L)^{O(1)}$
representing the coefficients of the polynomial $q_l$ with respect to $X_{r+1}$. 
Doing this for each $l$, $1\le l \le s$, again we have to add to the extended 
arithmetic network ${\cal N}_i$ some extra nodes which do not change its 
asymptotic size. Then we may suppose that ${\cal N}_i$ produces a division-free
straight-line program in $\mQ [X_1,\ldots , X_r]$ representing the coefficients
of the polynomials $q_1, \ldots q_s$ with respect to the variable $X_r$. As in the 
case of $i =0$ we discard by a call to the second subroutine the polynomials
$q_{m+1}, \ldots , q_s$ which do not have any zero in $\R^n$. 
>From the remaining polynomials $q_1, \ldots ,q_s$ we generate $q=q_1\cdots q_s$.
The additional costs of discarding $q_{m+1},\ldots , q_s$ and producing $q$ 
is of order $\left( \sum^s_{l=1}i d deg C_l L\right)^{O(1)} =
(i d deg Z L)^{O(1)}= (i d \delta^\ast_{i-1} L)^{O(1)}$. Thus,
without loss of generality we may suppose that the extended arithmetic network 
${\cal N}_i$ produces a division-free straight-line program in $\mQ [X_1, \ldots ,
X_r]$ of size  $\left( \sum^s_{l=1}i d deg C_l L \right)^{O(1)}=
 (i d \delta^\ast_{i-1} L)^{O(1)}$ which represents the coefficients of the 
polynomial $q$ with respect to the variable $X_r$. We observe that the point $t \in \mQ^r$
is a lifting point of the algebraic variety  $W^\ast_i = \cup^s_{l=1} C_l$, too.
Therefore, the lifting fiber of $t$ in $W^\ast_i$ is given by the rational number 
$\bar{\varrho}(t)$ and the coefficients of the polynomials $q(t,X_{r+1})$ and
$\bar{p}_1(t, X_{r+1}), \ldots, \bar{p}(t, X_{r+1})$, which, in principle, have
already been computed by the arithmetic network ${\cal N}_i$. Again applying
the procedure underlying \cite{gh3}, Theorem 28, to the input $\beta, t=(t_1,
\ldots , t_r), \bar{\varrho}(t), q(t,X_{r+1}), \bar{p}_1(t,X_{r+1}),
\ldots , \bar{p}_n(t, X_{r+1})$ we obtain a division-free straight-line program $\beta_i$
in $ \mQ [X_1, \ldots, X_r]$ of size $\Lambda_i = (i d \delta^\ast_i L)^{O(1)}$. The
straight-line program $\beta_i$ represents a non-zero polynomial $\varrho \in \mQ
[X_1, \ldots , X_r]$ and the coefficients with respect to $X_{r+1}$ of the 
polynomial $q$  and certain other polynomials $p_1, \ldots , p_n$ of $ \mQ [X_1,
\ldots , X_r, X_{r+1} ]$ having the properties (i) - (iv) stated in the 
Lemma 11.\\
The extended arithmetic network ${\cal N}_i$ over $\mQ$ which produces this 
output $\beta_i$ from the input $\beta_{i-1}$ and $\beta$ has size 
$(i d \delta^\ast_{i-1} L \Lambda_{i-1} )^{O(1)}$.\\
Observe that the length $\Lambda_i$ of the straight-line program $\beta_i$ is independent of the 
length $\Lambda_{i-1}$ of the input circuit $\beta_{i-1}$.  
More precisely, we have $\Lambda_i = (i d \delta_i^\ast L)^{O(1)}$. Taking into account
that $\delta^\ast_i \le d \delta^\ast_{i-1}$ and $\Lambda_{i-1} = ((i-1)d
\delta^\ast_{i-1}L)^{O(1)}$ holds we conclude 
that the size of the extended arithmetic network ${\cal N}_i$ which 
produces from the input circuits $\beta_{i-1}$ and $\beta$ the output circuits 
is of order $(i d \delta^\ast_i L)^{O(1)}$. Concatenating the networks
${\cal N}_1, \ldots , {\cal N}_i$ we finally obtain an extended arithmetic network 
${\cal N}$ over $\mQ$ which produces the straight-line program $\beta_i$ 
from the input circuit $\beta$. The network ${\cal N}$ is of size 
$(i d \delta^\ast L)^{O(1)}$.\hfill$\Box$\\

>From Lemma 11 one deduces now easily our main result.

\begin{theorem}\label{theorem12}
Let $n,d,\delta^\ast, L$ be natural numbers. Then there exists an extended
arithmetic network ${\cal N}$ over $\mQ$ of size $(n d \delta^\ast L)^{O(1)}$
with the following properties:\\
Let $f \in \mQ [X_1,\ldots ,X_n]$ be a non-constant polynomial of degree at 
most $d$ and suppose that $f$ is given by a division-free straight-line program
$\beta$ in $\mQ [ X_1, \ldots ,X_n]$ of length at most  $L$. Let 
$\Delta := \sum^n_{i=1} \left( \frac{\partial f}{\partial X_i} \right)^2$,
$W := \{ x \in \mC^n | f(x)=0\}$, $V:= W\cap \R^n = \{ x\in \R^n | f(x)=0\}$
and suppose that the variables $X_1, \ldots , X_n$ are in " sufficiently generic" position.
Furthermore, suppose that $V$ is a (non-empty) bounded smooth hypersurface of $\R^n$ with
regular equation $f$. For $0 \le i \le n$, let $W_i$ be the Zariski closure
of the set $ \{ x \in \mC^n | f(x) = \frac{\partial f (x)}{\partial X_1} =
\cdots = \frac{\partial f(x)}{\partial X_i} =0, \Delta (x) \not= 0\} $ and
$W^\ast_i := W_i \cap \R^n$.

Let $\delta^\ast_i := deg^\ast W_i := deg W^\ast_i$ be the real degree of the complex
variety  $W_i$ and assume that $\delta^\ast \ge \max \{ \delta^\ast_i |
0 \le i < n \}$ holds.\\
The algorithm represented by the extended arithmetic network ${\cal N}$ starts
from the straight-line program $\beta$ as input and produces a straight-line
program in $\mQ$ of size $(n d \delta^\ast L)^{O(1)}$. This straight-line program
represents the coefficients of $n+1$ univariate polynomials $q, p_1, \ldots , p_n \in \mQ
[X_n]$ satisfying the following conditions:
\begin{itemize}
\item[(i)] $deg q = \delta^\ast_{n-1}$
\item[(ii)] $\max \{ deg  p_i | 1 \le i \le n \} < \delta^\ast_{n-1}$
\item[(iii)] Any connected component of the real hypersurface $V$ has at least
one point contained in the set
\[ 
P := \{ (p_1(u), \ldots, p_n(u)) | u \in \R, q(u) =0\}.
\]
\end{itemize}
Moreover, the extended algorithmic network ${\cal N}$ codifies each real
zero $u$  of the polynomial $q$ (and hence, the elements of  $P$) " \`a la Thom".
\end{theorem}

{\bf Proof}:\\
Just apply Lemma 11 setting $i := n-1$. The remaining arguments are the 
same as in the proof of Theorem 8. \hfill $\Box$\\
 
\end{section}

\end{document}